\documentclass[12pt]{article}

\usepackage{amsmath, epsfig, cite}
\usepackage{amssymb}
\usepackage{amsfonts}
\usepackage{latexsym}
\usepackage{float}
\usepackage{mathrsfs}

\newtheorem{thm}{Theorem}[section]

\newtheorem{cor}[thm]{Corollary}

\numberwithin{equation}{section}

\UseRawInputEncoding
\begin{document}

%\linenumbers

\begin{center}
{\Large\bf On $\ell$-regular partitions and Hickerson's identity}
\end{center}

\vskip 2mm \centerline{Ji-Cai Liu}
\begin{center}
{\footnotesize Department of Mathematics, Wenzhou University, Wenzhou 325035, PR China\\
{\tt jcliu2016@gmail.com} }
\end{center}

%%date: November 27, 2014
%%date : December 4, 2014

\vskip 0.7cm \noindent{\bf Abstract.}
Based on two involutions and a bijection, we completely determine the difference between the number of $\ell$-regular partitions of $n$ into an even number of parts and into an odd number of parts for all positive integers $n$ and $\ell>1$, which extends two recent results due to Ballantine and Merca. As an application, we provide a combinatorial proof of Hickerson's identity on the number of partitions into an even and odd number of parts.

\vskip 2mm
\noindent{\it MR Subject Classifications}: 	05A17, 05A19
\section{Introduction}
A partition of a positive integer $n$ is a finite nonincreasing sequence of positive integers $(\lambda_1,\lambda_2,\cdots,\lambda_r)$ such that $\sum_{i=1}^r\lambda_i=n$.
For an integer $\ell>1$, a partition is called $\ell$-regular if none of its parts is divisible by $\ell$. Let $b_{\ell}(n)$ denote the number of the $\ell$-regular partitions of $n$. The generating function for $b_{\ell}(n)$ is given by
\begin{align*}
\sum_{n=0}^{\infty}b_{\ell}(n)q^n=\frac{(q^{\ell};q^{\ell})_{\infty}}{(q;q)_{\infty}},
\end{align*}
where $(a;q)_{\infty}=\prod_{k=0}^{\infty}(1-aq^k)$.
Let $b_{\ell}^e(n)$ be the number of $\ell$-regular partitions of $n$ into an even number of parts, and $b_{\ell}^o(n)$ be the number of $\ell$-regular partitions of $n$ into an odd number of parts.

Recently, Ballantine and Merca \cite{bm-qm-2023,bm-rasm-2023} showed that
\begin{align}
&b_{4}^e(n)-b_{4}^o(n)=(-1)^nd_{4}(n),\label{bm-1}\\[5pt]
&b_{6}^e(n)-b_{6}^o(n)=(-1)^nd_{6}(n),\label{bm-2}
\end{align}
where $d_{4}(n)$ is the number of partitions of $n$ into distinct parts which are not congruent to $2$ modulo $4$, and $d_{6}(n)$ is the number of partitions of $n$ into distinct parts which are not congruent to $\pm 2$ modulo $6$.

Inspired by Ballantine and Merca's combinatorial approach to \eqref{bm-1} and \eqref{bm-2}, we completely determine $b_{\ell}^e(n)-b_{\ell}^o(n)$ for all positive integers $n$ and $\ell>1$ through
two involutions and a bijection, which extends \eqref{bm-1} and \eqref{bm-2}.
\begin{thm}[Even case]\label{t-1}
Let $\ell>1$ be an even integer. For all positive integers $n$, we have
\begin{align}
b_{\ell}^e(n)-b_{\ell}^o(n)=(-1)^nd_{\ell}(n),\label{a-new-3}
\end{align}
where $d_{\ell}(n)$ denotes the number of partitions of $n$ into distinct parts which are congruent to $0,1,3,\cdots,\ell-3,\ell-1$ modulo $\ell$.
\end{thm}

\begin{thm}[Odd case]\label{t-2}
Let $\ell>1$ be an odd integer. For all positive integers $n$, we have
\begin{align}
b_{\ell}^e(n)-b_{\ell}^o(n)=(-1)^nc_{\ell}(n),\label{a-new-4}
\end{align}
where $c_{\ell}(n)$ denotes the number of partitions of $n$ into distinct odd parts which are not divisible by $\ell$.
\end{thm}

From \eqref{a-new-3} and \eqref{a-new-4}, we deduce the following result.
\begin{cor}
Let $\ell$ and $n$ be positive integers with $\ell>1$. For $\ell\equiv 0\pmod{2}$,
$b_{\ell}(n)$ and $d_{\ell}(n)$ have the same parity. For $\ell\equiv 1\pmod{2}$,
$b_{\ell}(n)$ and $c_{\ell}(n)$ have the same parity.
\end{cor}

For positive integers $n$ and $r$, let $q_r^e(n)$ be the number of partitions of $n$ into an even number of parts where each part occurs at most $r$ times, $q_r^o(n)$ be the number of partitions of $n$ into an odd number of parts where each part occurs at most $r$ times. Let $\Delta_r(n)=q_r^e(n)-q_r^o(n)$.
It was proved by Euler that
\begin{align}
\Delta_1(n)=
\begin{cases}
(-1)^j\quad&\text{if $n=j(3j\pm 1)/2$ for some $j=0,1,2,\cdots$},\\[5pt]
0\quad &\text{otherwise,}
\end{cases}\label{a-new-1}
\end{align}
which is known as the famous pentagonal number theorem (see \cite[page 10]{andrews-b-1998}).

By using the generating function technique, Hickerson \cite{hickerson-jcta-1973} showed that
\begin{align}
\Delta_3(n)=
\begin{cases}
(-1)^n\quad&\text{if $n=j(j+1)/2$ for some $j=0,1,2,\cdots$},\\[5pt]
0\quad &\text{otherwise,}
\end{cases}\label{a-new-2}
\end{align}
and established a simple formula for $\Delta_r(n)$ for $r\equiv 0\pmod{2}$ as follows.
\begin{thm}[Hickerson]\label{t-3}
For positive integers $n$ and positive even integers $r$, we have
\begin{align}
\Delta_r(n)=(-1)^nc_{r+1}(n).\label{a-1}
\end{align}
\end{thm}

For odd integers $r>3$, the formulae for $\Delta_r(n)$ have also been determined by several authors.
Alder and Muwafi \cite{am-fq-1795} established formulae for $\Delta_r(n)$ for $r=5,7$ in terms of the number of partitions into distinct parts taken from certain sets.
Hickerson\cite{hickerson-fq-1978} obtained a general formula for $\Delta_r(n)$ for any positive odd integer $r$. Robbins \cite{robbins-fq-2002} got a more simple formula for $\Delta_r(n)$ when $r$ is a positive odd integer.

We remark that the best known combinatorial approach to \eqref{a-new-1} is Franklin's involution proof of Euler's pentagonal number theorem (see \cite[page 10]{andrews-b-1998}). Fink, Guy and Krusemeyer \cite{fgk-cdm-2008} presented a combinatorial proof of \eqref{a-new-2} through a bijection of partitions
supplied with extra structure.

Notice that Hickerson's proof of \eqref{a-1} relies on the generating functions.
The second aim of the paper is to provide a combinatorial proof of \eqref{a-1} by using \eqref{a-new-4} and a bijection due to Glaisher \cite{lehmer-bams-1946}.
We shall present combinatorial proofs of \eqref{a-new-3}, \eqref{a-new-4} and \eqref{a-1} in Sections
2--4, respectively.

\section{Proof of Theorem \ref{t-1}}
For an $\ell$-regular partition $\lambda$ of $n$, let $e_{\ell}(\lambda)$ denote the largest even part of $\lambda$ which appears an odd number of times, and $t_{\ell}(\lambda)$ denote the largest repeated part of $\lambda$ which is not congruent to $\ell/2$ modulo $\ell$. For the sake of convenience, let $e_{\ell}(\lambda)=0$ for a partition $\lambda$ with all even parts occurring an even number of times, and $t_{\ell}(\lambda)=0$ for a partition $\lambda$ with no repeated part
$\not\equiv \ell/2\pmod{\ell}$.
Let $B'_{\ell}(n)$ be the set of $\ell$-regular partitions $\lambda$ of $n$ with $e_{\ell}(\lambda)\not=0$ or $t_{\ell}(\lambda)\not=0$. We define the involution $\psi_{\ell}$ on $B'_{\ell}(n)$ as follows.

{\bf \noindent Case 1.} If $2t_{\ell}(\lambda)>e_{\ell}(\lambda)$, then we merge two copies of $t_{\ell}(\lambda)$ into one part of double size. Since $t_{\ell}(\lambda)\not\equiv \ell/2\pmod{\ell}$, we have $\ell\nmid 2t_{\ell}(\lambda)$.\\
{\bf \noindent Case 2.} If $2t_{\ell}(\lambda)\le e_{\ell}(\lambda)$, then we split the part $e_{\ell}(\lambda)$ into two copies of $e_{\ell}(\lambda)/2$. Since $\ell\nmid e_{\ell}(\lambda)$, we have $\ell \nmid e_{\ell}(\lambda)/2$.

Let $l(\lambda)$ denote the number of the parts of the partition $\lambda$ (the length of $\lambda$).
It is clear that the involution $\psi_{\ell}$ reverses the parity of the length of the partitions.
For example,
\begin{align*}
&\psi_6: \quad (5,4,3^2,2^2,1)\longleftrightarrow (5,3^2,2^4,1),\\[5pt]
&\psi_8: \quad (7,5,4^3,2^3,1) \longleftrightarrow (7,5,4^2,2^5,1).
\end{align*}

Let $A_{\ell}(n)$ be the set of partitions of $n$ into distinct odd parts and the parts $\equiv \ell/2\pmod{\ell}$ with even multiplicity for $\ell \equiv 0\pmod{4}$, and the set of partitions of $n$ into distinct odd parts $\not\equiv \ell/2\pmod{\ell}$ and the parts $\equiv \ell/2\pmod{\ell}$ with arbitrary multiplicity for $\ell \equiv 2\pmod{4}$.
Notice that $e_{\ell}(\lambda)=t_{\ell}(\lambda)=0$ if and only if $\lambda\in A_{\ell}(n)$.

It is clear that $l(\lambda)\equiv n\pmod{2}$ for $\lambda\in A_{\ell}(n)$.
By this fact and the above involution $\psi_{\ell}$ on $B'_{\ell}(n)$, to prove \eqref{a-new-3}, it suffices to show that $\# A_{\ell}(n)=\#D_{\ell}(n)$, where $D_{\ell}(n)$ denotes the set of partitions of $n$ into distinct parts which are congruent to $0,1,3,\cdots,\ell-3,\ell-1$ modulo $\ell$. We shall establish a bijection $\sigma_{\ell}: A_{\ell}(n)\to D_{\ell}(n)$.

{\noindent \it The map $\sigma_{\ell}$ from $A_{\ell}(n)$ to $D_{\ell}(n)$:} If the partition contains two copies of the same part, then we merge the two parts into one part of double size. We repeat this procedure until all parts are distinct.
For $\ell \equiv 0\pmod{4}$ and $\lambda\in A_{\ell}(n)$, the even parts $\equiv \ell/2\pmod{\ell}$ of $\lambda$ appear an even number of times.
For $\ell \equiv 2\pmod{4}$ and $\lambda\in A_{\ell}(n)$, the odd parts $\equiv \ell/2\pmod{\ell}$ appear an arbitrary number of times. It follows that $\sigma_{\ell}\left(\lambda\right)\in D_{\ell}(n)$ for
$\lambda\in A_{\ell}(n)$.

For example,
\begin{align*}
&\sigma_6:\quad (9^2,7,5,3^3,1)\to (18,7,6,5,3,1),\\[5pt]
&\sigma_8:\quad (12^2,9,7,5,4^4,3)\to (24,9,8^2,7,5,3)\to (24,16,9,7,5,3).
\end{align*}

{\noindent \it The map $\sigma_{\ell}^{-1}$ from $D_{\ell}(n)$ to $A_{\ell}(n)$:}
If the partition contains a part $\equiv 0\pmod{\ell}$, then we split this part into two equal halves. We repeat this procedure until no part is congruent to $0$ modulo $\ell$.
For $\ell \equiv 0\pmod{4}$ and $\lambda\in D_{\ell}(n)$, the even parts $\equiv \ell/2\pmod{\ell}$ do not appear in $\lambda$.
For $\ell \equiv 2\pmod{4}$ and $\lambda\in D_{\ell}(n)$, the odd parts $\equiv \ell/2\pmod{\ell}$ appear
at most once in $\lambda$. It follows that $\sigma_{\ell}^{-1}\left(\lambda\right)\in A_{\ell}(n)$ for
$\lambda\in D_{\ell}(n)$.

For example,
\begin{align*}
&\sigma_6^{-1}:\quad (18,7,6,5,3,1)\to (9^2,7,5,3^3,1),\\[5pt]
&\sigma_8^{-1}:\quad (24,16,9,7,5,3)\to (12^2,9,8^2,7,5,3)\to (12^2,9,7,5,4^4,3).
\end{align*}

\section{Proof of Theorem \ref{t-2}}
For an $\ell$-regular partition $\lambda$ of $n$, let $e_{\ell}(\lambda)$ denote the largest even part of $\lambda$ which appears only once, and $t_{\ell}(\lambda)$ denote the largest repeated part of $\lambda$.
For the sake of convenience, let $e_{\ell}(\lambda)=0$ for a partition $\lambda$ with no even part and $t_{\ell}(\lambda)=0$ for a partition $\lambda$ with distinct parts.
Let $B'_{\ell}(n)$ be the set of $\ell$-regular partitions $\lambda$ of $n$ with $e_{\ell}(\lambda)\not=0$ or $t_{\ell}(\lambda)\not=0$. We define the involution $\psi_{\ell}$ on $B'_{\ell}(n)$ as follows.

{\noindent \bf Case 1.} If $2t_{\ell}(\lambda)>e_{\ell}(\lambda)$, then we merge two copies of $t_{\ell}(\lambda)$ into one part of double size. Since $t_{\ell}(\lambda)$ is not divisible by the odd integer $\ell$, we have $\ell\nmid 2t_{\ell}(\lambda)$.\\
{\noindent \bf Case 2.} If $2t_{\ell}(\lambda)\le e_{\ell}(\lambda)$, then we split the part $e_{\ell}(\lambda)$ into two copies of $e_{\ell}(\lambda)/2$. Since $e_{\ell}(\lambda)$ is not divisible by the odd integer $\ell$, we have $\ell\nmid e(\lambda)/2$.

It is clear that the involution $\psi_{\ell}$ reverses the parity of the length of the partitions.
We illustrate this involution through the following example $n=10$ and $\ell=3$:
\begin{align*}
&(4,1^6)\longleftrightarrow (2^2,1^6)& &(4,2^2,1^2) \longleftrightarrow (2^4,1^2)\\[5pt]
&(4^2,2) \longleftrightarrow (8,2)&  &(5,2,1^3)\longleftrightarrow (5,1^5)\\[5pt]
&(5,4,1)\longleftrightarrow (5,2^2,1)& & (7,2,1)\longleftrightarrow (7,1^3)\\[5pt]
&(8,1^2) \longleftrightarrow (4^2,1^2)& & (10) \longleftrightarrow (5^2)\\[5pt]
&(2,1^8) \longleftrightarrow (1^{10})&& (2^3,1^4) \longleftrightarrow (4,2,1^4)\\[5pt]
&(2^5) \longleftrightarrow (4,2^3)&&
\end{align*}

Let $A_{\ell}(n)$ be the set of partitions of $n$ into distinct odd parts which are not divisible by $\ell$. Note that $e_{\ell}(\lambda)=t_{\ell}(\lambda)=0$ if and only if $\lambda\in A_{\ell}(n)$. It is clear that $n\equiv l(\lambda)\pmod{2}$ for $\lambda\in A_{\ell}(n)$. Then the proof of \eqref{a-new-4} follows from this fact and the above involution $\psi_{\ell}$ on $B'_{\ell}(n)$.

\section{Proof of Theorem \ref{t-3}}
Let $Q_r(n)$ be the set of partitions of $n$ where each part occurs at most $r$ times,
and $B_{\ell}(n)$ be the set of $\ell$-regular partitions of $n$.
There is a classic bijection between $Q_r(n)$ and $B_{r+1}(n)$ due to Glaisher \cite{lehmer-bams-1946}, which is stated as follows.

{\noindent\it The map $\varphi_r$ from $Q_r(n)$ to $B_{r+1}(n)$:}
If the partition contains a part which is a multiple of $r+1$, then we split this part into $r+1$ equal parts. We repeat this procedure until no part is a multiple of $r+1$.

{\noindent\it The map $\varphi_r^{-1}$ from $B_{r+1}(n)$ to $Q_r(n)$:}
If the partition contains $r+1$ copies of the same part, then we merge the $r+1$ parts into one part of $(r+1)$-fold size. We repeat this procedure until no part appears more than $r$ times.

Note that each procedure in $\varphi_r$ splits one part into an odd number of parts, which preserves the parity of the length of the partitions. It follows that $l\left(\varphi_r(\lambda)\right)\equiv l(\lambda)\pmod{2}$ for $\lambda\in Q_r(n)$.

For example,
\begin{align*}
\varphi_2:&\quad (9,6^2,3,2^2,1)\to (3^3,2^8,1^4)\to (2^8,1^{13}),\\[5pt]
\varphi_2^{-1}:&\quad (2^8,1^{13})\to (6^2,3^4,2^2,1)\to (9,6^2,3,2^2,1).
\end{align*}

The proof of \eqref{a-1} follows from \eqref{a-new-4} and the above bijection between $Q_r(n)$ and $B_{r+1}(n)$ which preserves the parity of the length of the partitions.

\vskip 5mm \noindent{\bf Acknowledgments.} This work was supported by the National Natural Science Foundation of China (grant 12171370).

\end{document}